\documentclass[12pt]{ttuthes}

\usepackage{amsmath,amssymb,latexsym,amsthm}
\usepackage{graphicx}

\newtheorem{theorem}{Theorem}[section]

\renewcommand{\thechapter}{\Roman{chapter}}

\begin{document}

\frontmatter

\begin{center}AREA OF POLYGONS IN HYPERBOLIC GEOMETRY\\
By\\
BENJAMIN AARON BAILEY\\
A THESIS\\
IN\\
MATHEMATICS\\
Submitted to the Graduate Faculty\\
of Texas Tech University in \\
Partial Fulfillment of\\
the Requirements for\\ 
the Degree of\\
\medskip
MASTER OF SCIENCE
\end{center}

\chapter{ACKNOWLEDGEMENTS}
My sincerest thanks and appreciation goes to my thesis committee: Dr.\ Roger Barnard, Dr.\ Petros Hadjicostas (chair),
 Dr.\ Magdalena Toda, and Dr.\ Brock Williams for their support and knowledge offered over the last semester.

\tableofcontents

\chapter{ABSTRACT}
\noindent  Consider the Poincare model for hyperbolic geometry on the unit disc, and an arbitrary $n$-gon in this 
geometry.  Chapter I gives a briefly introduces the conformal metric which generates hyperbolic geometry.   
In Chapter II, there is an analytic derivation of a convenient computational formula for the hyperbolic area 
in terms of the coordinates of its vertices, as well as an insightful geometric interpretation of this formula in terms of 
naturally occuring angles of the $n$-gon.  Chapter III extends the formulas in Chapter II to the closed unit disc, with an 
alternative, more geometrically motivated proof.  Chapter IV uses the results of Chapter I to establish identities between 
hyperbolic area and perimeter of an $n$-gon.  A proof of the existence of a solution of the following isoperimetric problem 
is also given:  Maximize (if such a maximum exists) the area of an $n$-gon with fixed perimeter.
            
\listoffigures
\addcontentsline{toc}{chapter}{LIST OF FIGURES}

\mainmatter

\setcounter{page}{1}
\pagenumbering{arabic}

\chapter{INTRODUCTION}
\setcounter{equation}{0}

\noindent Consider the Poincare model of hyperbolic geometry on the unit disc $\{z \in \mathbb{C} :|z| \leq 1\}$.
  According to \cite{Anderson}, p.97, the hyperbolic arclength $P$ of a curve $\gamma$ of bounded variation  is 
defined by   
\begin{equation}\label{perimeter}
P=\int\limits_\gamma \frac{2}{1-|z|^2}|\mathrm{d}z|,
\end{equation}
Also, by  \cite{Anderson}, p.138, the hyperbolic area of a region $\mathcal{R}$ is given by
\begin{equation}\label{area}
A=\iint\limits_\mathcal{R}\left(\frac{2}{1-|z|^2}\right)^2  \mathrm{d}x \mathrm{d}y=
 \iint\limits_\mathcal{R} \left(\frac{2}{1-x^2-y^2}\right)^2\, \mathrm{d}x \mathrm{d}y.
\end{equation}
Equation~(\ref{perimeter}) in differential form is
\begin{equation}
\mathrm{d}P = \frac{2}{1-|z|^2}|\mathrm{d}z| = \rho(z)|\mathrm{d}z|.\nonumber
\end{equation}

\noindent This gives a function $\rho(z)$ which relates the hyperbolic separation $\mathrm{d}P$ of two nearby points 
to their Euclidean separation $|\mathrm{d}z|.$  In a certain sense, this function describes the ``density" of space at a 
given point.  This function serves to define distance and the notion of a line (more properly called a geodesic) in this 
geometry.  Take $z$ and $w$ in the interior of the unit disc.  By  \cite{Anderson}, p.89, the distance between these two 
points is the infimum of the lengths (according to equation~(\ref{perimeter})) of all curves joining these two points.  
A geodesic between $z$ and $w$ is a curve whose length attains this infimum.   From \cite{Needham} p.316, there is a 
unique geodesic between two points, and this geodesic is an arc of a circle which intersects the unit circle 
at right angles.  If one of the points lies on the boundary, the geodesic can still be defined in this geometric sense 
even though the distance to any other point is infinite.\\
We can say that $\rho(z)=\frac{2}{1-|z|^2}$ generates hyperbolic geometry as it completely determines geodesics 
between all points in the unit disc.  This background prepares us for the following formulas and results. 

\chapter[ANALYTICAL DERIVATION OF MAIN RESULTS]{ANALYTICAL DERIVATION AND GEOMETRIC INTERPRETATION 
OF A COMPUTATIONAL FORMULA FOR HYPERBOLIC AREA OF $N$-GONS}

\setcounter{equation}{0}
\setcounter{theorem}{0}

\noindent From \cite{Anderson}, p.137, the formula for area of a hyperbolic $n$-gon (henceforth referred to as an $n$-gon) 
is given by 
\begin{equation}\label{classical}
A=\pi(n-2)-\sum_{k=1}^n \theta_k,
\end{equation}
 where $\theta_k$ represents the interior angle of the $k^{th}$ vertex.  This formula, though geometrically elegant, is 
not computational in nature unless the interior angles are given.  The following theorem provides a simple 
computational formula for the hyperbolic area (henceforth referred to as area), where the coordinates of the vertices are
 given, that also yields an equally elegant geometric interpretation.  The following proof and interpretation assume that
 the vertices are in the interior of the unit disc, though in Chapter III a proof will be given that removes this restriction.
\begin{theorem}
Consider the Poincare disc model for hyperbolic geometry.  Let $(x_1,y_1), \ldots , (x_n,y_n)$ defined in the interior of 
 the unit disc be the positively oriented vertices of  an n-gon.  Define $(x_{n+1},y_{n+1})$ to be $(x_1,y_1)$.  The area is 
given by the formula
\end{theorem}
\begin{equation}\label{compform}
A =\sum_{k=1}^n 2\;{\rm Tan}^{-1} \left( {\frac{x_k y_{k+1}-y_k x_{k+1}}{1-x_k x_{k+1}-y_k y_{k+1}}} \right).
\end{equation}
\noindent \textit{Proof.} We seek to evaluate equation~(\ref{area}) where $\mathcal{R}$ is the region bounded by 
$\gamma$, the boundary of the $n$-gon.
\noindent The proof is divided as follows:\\
1.  Conversion of the area double integral to a line integral via Green's theorem.\\
2.  Determining a parameterization $\gamma$ for the boundary of an $n$-gon.\\
3.  Evaluation of the area as a line integral over $\gamma$.\\
\noindent 1.)  Green's theorem in the plane states that if  $\textbf{B} = M(x,y) \textbf{i}+N(x,y) \textbf{j}$ is a
 continuously differentiable vector field on $\mathcal{R}$  and $\gamma \equiv \partial \mathcal{R}$, the boundary of 
$\mathcal{R}$, is  a Jordan curve of bounded variation, then
\begin{equation}\label{green}
\iint\limits_\mathcal{R} \left(\frac{\partial N}{\partial x}-\frac{\partial M}{\partial y}\right) \mathrm {d}x \mathrm {d}y=
\oint\limits_{\partial R} M \mathrm {d}x+N\mathrm {d}y.
\end{equation}
\noindent Applying Green's formula for $x^2+y^2 < 1$, let
\begin{equation}
\textbf{B}=M(x,y) \textbf{i}+N(x,y) \textbf{j}=\frac{-2y}{1-x^2-y^2}\textbf{i}+\frac{2x}{1-x^2-y^2}\textbf{j}.
\end{equation}
\noindent Substituting this quantity into equation~(\ref{green}) we obtain
\begin{equation}\label{lineint1}
A=\iint\limits_\mathcal{R}
\left(\frac{2}{1-x^2-y^2}\right)^2 \mathrm {d}x \mathrm {d}y = \oint\limits_\gamma \frac{2}{1-x^2-y^2}
(x\mathrm {d}y-y\mathrm{d}x).
\end{equation}
\noindent  Consider $\gamma \equiv \partial \mathcal{R}$ as a subset of $\mathbb{C}$, then $x^2+y^2=|z|^2$
 and $\bar{z}\mathrm {d}z = (x-iy)(\mathrm{d}x+i\mathrm {d}y)$ implies
$x \mathrm {d}y-y \mathrm {d}x ={\rm Im} (\bar{z} \mathrm {d}z)$ so that the right hand side of the area 
equation~(\ref{lineint1}) becomes
\begin{equation}\label{lineint2}
A={\rm Im}\oint\limits_\gamma \frac{2\bar{z}}{1-|z|^2} \mathrm {d}z.
\end{equation}
\noindent The above expression for area is sufficient for the proof of equation~(\ref{compform}), but a stronger statement is 
more fruitful for the later development of identities.\\
\noindent \textit{Claim:}  $I=\oint\limits_\gamma \frac{2\bar{z}}{1-|z|^2} \mathrm {d}z$ 
is purely imaginary.\\

\noindent \textit{Proof of claim:}  Clearly $\bar{I}=\oint\limits_\gamma \frac{2z}{1-|z|^2} \mathrm {d}\bar{z}$.  Adding and 
converting to real form, we have
\begin{equation}
I+\bar{I}=\oint\limits_\gamma \left(\frac{4x}{1-x^2-y^2}\mathrm{d}x+\frac{4y}{1-x^2-y^2}\mathrm{d}y\right).
\end{equation}
Applying Green's formula for $x^2+y^2 < 1$, let
\begin{equation}
\textbf{B}=M(x,y)\textbf{i}+N(x,y)\textbf{j}=\frac{4x}{1-x^2-y^2}\textbf{i}+\frac{4y}{1-x^2-y^2}\textbf{j}.
\end{equation}
Substituting this quantity into equation~(\ref{green}) we obtain
\begin{equation}
0=\oint\limits_\gamma \left(\frac{4x}{1-x^2-y^2}\mathrm{d}x+\frac{4y}{1-x^2-y^2}\mathrm{d}y\right)=I+\bar{I}.
\end{equation}

\noindent So we have ${\rm Re}(I)=(I+\bar{I})/2=0$ which implies $I=i\ {\rm Im}(I)$, that is $I$ is purely imaginary.  This 
concludes the proof of the claim.\\
The immediate consequence of this result is
\begin{equation}\label{pureimaginary}
A=\frac{1}{i}\oint\limits_\gamma \frac{2\bar{z}}{1-|z|^2} \mathrm {d}z.
\end{equation}
\noindent 2.)  By \cite{Stillwell}, p.89, the isometries in this model of hyperbolic geometry are the M\"{o}bius 
transformations of the following form,  $M_a(z)=e^{i\theta}\frac{z-a}{1-\bar{a}z}$ where $a$ is a point in the open unit 
disc and $\theta$ is a real number.  These mappings are disc automorphisms and form a group with respect to 
composition.  These properties will be used to parameterize the geodesic from 
$z_k=(x_k,y_k)$ to $z_{k+1}=(x_{k+1}, y_{k+1}).$  In this chapter it is assumed that consecutive points are distinct.\\
\noindent Let $\gamma_k:[0,1]\longrightarrow \mathbb{C}$ be the geodesic from $z_k$ to $z_{k+1}$ in the open unit disc.
  The hyperbolic arclength of $\gamma_k$, equation~(\ref{perimeter}), is invariant under the transformation 
$M_k(z)=\frac{z-z_k}{1-\bar{z}_kz}.$  Note that $M_k(z_k)=0$ and
 $M_k(z_{k+1})=\frac{z_{k+1}-z_k}{1-\bar{z}_kz_{k+1}} \neq 0$.  Define the quantity
\begin{equation}
a_k=\frac{1}{M_k(z_{k+1})}=\frac{1-\bar{z}_kz_{k+1}}{z_{k+1}-z_{k}}.  
\end{equation}
Now $|z_{k+1}|<1$ implies $|M_k(z_{k+1})|<1$ which forces $|a_k|>1.$
  Since $M_k$ is an isometry, the image of $\gamma_k$ under $M_k$ must be the geodesic from $0$ to 
$M_k(z_{k+1})$.  Geodesics from the origin to a point are Euclidean lines, so an appropriate equation for this 
ray is $M_k(\gamma_k(t))=t/a_k$ where $t \in [0,1]$.  We now have $\gamma_k(t)=M_k^{-1}(t/a_k)$, where 
$M_k^{-1}(z)=\frac{z+z_k}{z\bar{z}_k+1}$, so that we have the following parameterization
\begin{equation}\label{geodesic}
\gamma_k(t)=\frac{t+a_kz_k}{t\bar{z}_k+a_k} \qquad\mbox{where}\qquad t\in [0,1].
\end{equation}
\noindent 3.)  Given equation~(\ref{geodesic}) for each side of the $n$-gon, the right hand side of
equation~(\ref{pureimaginary}) can be expressed as
\begin{equation}
A=\frac{1}{i}\sum\limits_{k=1}^n \int\limits_{\gamma_k} \frac{2\bar{z}}{1-|z|^2}\mathrm {d}z,
\end{equation}
which after substitution becomes
\begin{equation}\label{ugly}
A=\frac{1}{i}\sum\limits_{k=1}^n \int_0^1 \frac{2\left(\frac{t+\bar{a}_k\bar{z}_k}{tz_k+\bar{a}_k}\right)\frac{{d}}{{d}t}
\left(\frac{t+a_kz_k}{t\bar{z}_k+a_k}\right)}{1-\left|\frac{t+a_kz_k}{t\bar{z}_k+a_k}\right|^2}\mathrm {d}t.
\end{equation}
Algebraic simplification of the integrand yields
\begin{equation}
A=\frac{1}{i}\sum\limits_{k=1}^n \int_0^1 \frac{2a_k(t+\bar{a}_k\bar{z}_k)}{(t\bar{z}_k+a_k)(|a_k|^2-t^2)}\mathrm {d}t.
\end{equation}
After partial fraction decomposition we have
\begin{equation}\label{integralsum}
A=\frac{1}{i}\sum\limits_{k=1}^n \int_0^1 \left[\frac{2\bar{z}_k/a_k}{t\bar{z}_k/a_k+1}-\frac{1}{|a_k|+t}+\frac{1}{|a_k|-t}\right]\mathrm {d}t.
\end{equation}
For the first integrand recall that $|a_k|>1$ so that for all $t \in [0,1]$ 
we have $\left|\frac{\bar{z}_kt}{a_k}\right|<1.$  This shows that the denominator has a positive 
real part, so that we may integrate using the principal branch of the logarithm.  Combining this with the last two 
integrands we have
\begin{equation}\label{areacomplex}
A=\frac{1}{i}\sum\limits_{k=1}^n 2\log\left(1+\frac{\bar{z}_k}{a_k}\right)+
\frac{1}{i}\sum\limits_{k=1}^n \ln\left(\frac{|a_k|^2}{|a_k|^2-1}\right).
\end{equation}
Gathering real and imaginary parts,
\begin{equation}\label{areacomplex2}
A=\sum\limits_{k=1}^n 2\arg\left(1+\frac{\bar{z}_k}{a_k}\right)+
\frac{1}{i}\left[\sum\limits_{k=1}^n \ln\left|1+\frac{\bar{z}_k}{a_k}\right|^2+
\sum\limits_{k=1}^n \ln\left(\frac{|a_k|^2}{|a_k|^2-1}\right)\right]
\end{equation}
which implies
\begin{eqnarray}
A &=&\sum\limits_{k=1}^n 2\arg\left(1+\frac{\bar{z}_k}{a_k}\right)\label{argformula} = 
\sum\limits_{k=1}^n 2\arg\left(\frac{1-|z_k|^2}{1-\bar{z}_kz_{k+1}}\right)\nonumber\\ &=& 
\sum\limits_{k=1}^n 2\arg\left(\frac{1-|z_k|^2}{|1-\bar{z}_kz_{k+1}|^2}(1-z_k\bar{z}_{k+1})\right)
=\sum\limits_{k=1}^n 2\arg\left(1-z_k\bar{z}_{k+1}\right)\nonumber\\ &=&
\sum\limits_{k=1}^n 2\ {\rm Tan}^{-1}\left( \frac{{\rm Im}(1-z_k\bar{z}_{k+1})}{{\rm Re}(1-z_k\bar{z}_{k+1})} \right)
\quad\mbox{(since ${\rm Re} (1-z_k\bar{z}_{k+1})>0$);}\nonumber\\
A &=&\sum_{k=1}^n 2\;{\rm Tan}^{-1} \left( {\frac{x_k y_{k+1}-y_k x_{k+1}}{1-x_k x_{k+1}-y_k y_{k+1}}} \right).\nonumber
\end{eqnarray}
$\Box$

\noindent We now have a convenient computational formula for area, however, the use of
$\arg$ and ${\rm Tan}^{-1}$ suggests the geometric interpretation should be based on naturally occuring 
angles in the $n$-gon. This is the case, but before the formula is given, some definitions
need to be introduced.\\
\noindent \textit{Definition 1.)}  A side of an $n$-gon is said to be an \textit{inward side} if it is concave with respect to the 
interior of the $n$-gon in the Euclidean sense.\\
\noindent \textit{Definition 2.)}  A side of an $n$-gon is said to be an \textit{outward side} if it is convex with respect to the 
interior of the $n$-gon in the Euclidean sense.  A side that corresponds to a Euclidean line segment will be considered 
an outward side.\\

\begin{figure}
\setlength{\unitlength}{1cm}
\begin{picture}(6.0,6.0)(-5,-5)
\rotatebox{270}{\includegraphics{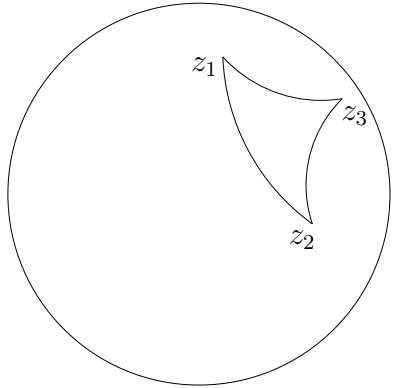}}
\end{picture}
\put(1.1,4.65){$z_1$}
\put(2.4,2.35){$z_2$}
\put(3.1,4){$z_3$}
\caption{Illustration of inward and outward sides.}\label{sidedef}
\end{figure}

\begin{figure}
\setlength{\unitlength}{1cm}
\begin{picture}(6.0,6.0)(-5,-5)
\rotatebox{270}{\includegraphics{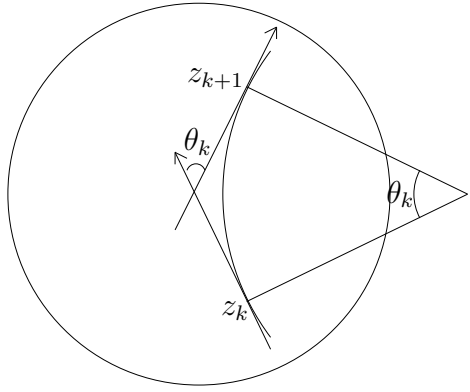}}
\end{picture}
\put(3.7,2.9){$\theta_k$}
\put(1,3.6){$\theta_k$}
\put(1.5,1.4){$z_k$}
\put(1.05,4.5){$z_{k+1}$}
\caption{Angle $\theta_k$ traversed by the tangent vector.}\label{traversed}
\end{figure}

\noindent For example, in Figure ~\ref{sidedef}, the geodesic from $z_1$ to $z_2$ is outward, and the geodesic from 
$z_2$ to $z_3$ is inward.
\begin{theorem}
Consider an $n$-gon with $l$ inward sides and $m$ outward sides.  Then the
area of the $n$-gon is given by the formula
\begin{equation}\label{geometric}
A=\sum_{k=1}^l \alpha_k-\sum_{k=1}^m \beta_k
\end{equation}
where $\alpha_k$ and $\beta_k$ denote the angle measure of the arcs of the inward and outward sides, respectively.
\end{theorem}
\noindent The explanation of this theorem calls for an analytic condition for a side to be inward or outward.  Refer to 
Figure ~\ref{sidedef}.  A side is inward if the smaller angle (magnitude-wise, measured in polar coordinates and 
from the origin) swept from $z_k$ to $z_{k+1}$ is positive but not equal to $\pi$.  Similarly, a side is outward if the 
smaller angle (magnitude-wise, measured in polar coordinates and from the origin) swept from $z_k$ to $z_{k+1}$ is 
negative, $0$ or $\pm\pi$.  Let $\hat{\textbf{k}}$ be the positive normal to the complex plane.  If $z_k$, $z_{k+1}$ are 
visualized as vectors in the plane, then the sign of the cross product 
$(z_k \times z_{k+1})\cdot\hat{\textbf{k}}=x_ky_{k+1}-y_kx_{k+1}$ in the 
$\hat{\textbf{k}}$ direction will determine the side's inwardness (sign is positive) or outwardness (sign is not positive), 
so that 
\begin{eqnarray}
x_ky_{k+1}-y_kx_{k+1} & > &0 \quad\mbox{implies}  \quad \gamma_k \quad\mbox{is inward, and}\label{criterion1}\\
x_ky_{k+1}-y_kx_{k+1} & \leq & 0 \quad\mbox{implies} \quad  \gamma_k \quad\mbox{is outward}\label{criterion2}.
\end{eqnarray}
Note that this condition is valid even if the vertices lie on the boundary. This implies that a term in the sum of 
equation~(\ref{compform}) is positive if $\gamma_k$ is inward, and non-positive if $\gamma_k$ is outward.  The sum 
may be written as
\begin{equation}\
A =\sum_{\gamma_k -in} \left|2\;{\rm Tan}^{-1} \left( {\frac{x_k y_{k+1}-y_k x_{k+1}}{1-x_k x_{k+1}-y_k y_{k+1}}} \right)\right|
-\sum_{\gamma_k -out} \left|2\;{\rm Tan}^{-1} \left( {\frac{x_k y_{k+1}-y_k x_{k+1}}{1-x_k x_{k+1}-y_k y_{k+1}}} \right)\right|,
\end{equation}
or equivalently as
\begin{equation}\label{geometricanalytic}
A=\sum_{\gamma_k -in} \left|2\arg\left(1+\frac{\bar{z}_k}{a_k}\right)\right|-
\sum_{\gamma_k -out} \left|2\arg\left(1+\frac{\bar{z}_k}{a_k}\right)\right|.
\end{equation}
The only step that remains is to geometrically interpret the angles $\left|2\arg\left(1+\frac{\bar{z}_k}{a_k}\right)\right|$.  
Since the geodesics are actually segments of circles, the angular measure of the arc from $z_k$ to $z_{k+1}$
is simply the angle swept out by the tangent vector from $t=0$ to $t=1$, that is the magnitude of the angle between the 
vectors $\gamma_k^\prime(0)$ and $\gamma_k^\prime(1)$ (Figure~\ref{traversed}).

\noindent Represent these vectors by $(a_1,b_1)$ and $(a_2,b_2)$ 
respectively, then the cosine of the angle $\theta_k$ between them is 
$\cos(\theta_k)=\frac{a_1a_2+b_1b_2}{\sqrt{a_1^2+b_1^2}\sqrt{a_2^2+b_2^2}}$.  If the vectors are represented by 
$a_1+ib_1$ and $a_2+ib_2$ then
\begin{equation}
\cos(\theta_k)=
\frac{{\rm Re}(\gamma_k^\prime(0)\overline{\gamma_k^\prime(1)})}{\left|\gamma_k^\prime(0)\right|\left|
\gamma_k^\prime(1)\right|}.
\end{equation}
Substituting the values for $\gamma_k^\prime(0)$ and $\gamma_k^\prime(1)$ we obtain
\begin{eqnarray}
\cos(\theta_k) &=& {\rm Re}\left(\frac{\left(\frac{1-|z_k|^2}{a_k}\right)\left(\frac{\bar{a}_k(1-|z_k|^2)}{(z_k+\bar{a}_k)^2}\right)}{\left|
\frac{1-|z_k|^2}{a_k}\right|\left|\frac{a_k(1-|z_k|^2)}{(\bar{z}_k+a_k)^2}\right|}\right)
= {\rm Re}\left(\frac{\bar{a}_k|\bar{z}_k+a_k|^2}{a_k(z_k+\bar{a}_k)^2}\right)\nonumber\\
&=& {\rm Re}\left(\frac{|a_k|^2\bar{a}_k\left|1+\frac{z_k}{\bar{a}_k}\right|^2}{a_k\bar{a}_k^2
\left(1+\frac{z_k}{\bar{a}_k}\right)^2}\right)
= {\rm Re}\left(\frac{\left|1+\frac{\bar{z}_k}{a_k}\right|^2\left(1+\frac{\bar{z}_k}{a_k}\right)^2}{\left(
1+\frac{z_k}{\bar{a}_k}\right)^2\left(1+\frac{\bar{z}_k}{a_k}\right)^2}\right)\nonumber\\
&=& {\rm Re}\left(\frac{\left|1+\frac{\bar{z}_k}{a_k}\right|^2}{\left|1+\frac{\bar{z}_k}{a_k}\right|^4}
\left(1+\frac{\bar{z}_k}{a_k}\right)^2\right) = 
{\rm Re} \left(\frac{\left(1+\frac{\bar{z}_k}{a_k}\right)^2}{\left|1+\frac{\bar{z}_k}{a_k}\right|^2}\right)\nonumber\\
&=& \cos\left(\arg\left(1+\frac{\bar{z}_k}{a_k}\right)^2\right) = \cos\left|2\arg\left(1+\frac{\bar{z}_k}{a_k}\right)\right|.
\end{eqnarray}
Now both $\theta_k$ and $\left|2\arg\left(1+\frac{\bar{z}_k}{a_k}\right)\right|$ are between $0$ and $\pi$, so we conclude 
that
\begin{equation}
\theta_k = \left|2\arg\left(1+\frac{\bar{z}_k}{a_k}\right)\right|.
\end{equation}
By virtue of equation~(\ref{geometricanalytic}), the theorem is proven, so that we have an expression 
for area with a convenient analytic and geometric form. $\Box$

\chapter[GEOMETRIC DERIVATION OF THE MAIN RESULTS]{AN ALTERNATIVE DERIVATION AND EXTENSION OF THE 
COMPUTATIONAL AND GEOMETRIC FORMULAS}
\setcounter{equation}{0}
\noindent  The derivation in Chapter II is based on an analytic approach toward hyperbolic geometry.  The alternative 
derivation of equations~(\ref{geometric}) and~(\ref{compform}) in this chapter is more geometrically motivated, and 
results in an extension of these equations to the closed unit disc.  The respective 
strengths and weaknesses of these methods of proof will be discussed in the concluding remarks.  The following theorem
 (modified for our purposes) is found in \cite{Do Carmo}, p.250.
\begin{theorem}\label{theta}
Let $z(t)=a(t)+ib(t)$, $t\in [t_0,t_1]$, with $a,b \in C^1[t_0,t_1]$ and $|z(t)|$ identically $1$.  There exists $\phi(t) \in
 C^1[t_0,t_1]$ with $z(t)= e^{i\phi(t)}.$ 
\end{theorem}
\noindent  Note that the difference between any two functions $\phi(t)$ is a constant multiple of $2\pi$.  In the context of 
our problem $z(t)$ will be the tangent vector $\textbf{T}_k(t) = \frac{\gamma_k^\prime(t)}{|\gamma_k^\prime(t)|}.$  We 
have only defined a parameterization when all vertices are in the interior of the unit disc, but certainly, a 
parameterization exists when this restriction is removed.  Essentially, $\phi(t)$ is the differentiable function yielding 
the angle with respect to the x-axis of $\textbf{T}_k(t)$ as it travels over the boundary of the $n$-gon.  The boundary is 
piecewise $C^2$ so that except at vertices, $\phi(t)$ exists.  Adhering to convention, clockwise angles will be negative, 
and counterclockwise angles will be positive.  Define exterior angles at the vertices in the following way:
\begin{equation}\label{thetaext}
\theta_{ext} = \pi -\theta_{int},
\end{equation}
where $\theta_{int}$ is the positive interior angle at each vertex.  These quantities are also defined for vertices which lie 
on the boundary of the unit disc.  Note that $\theta_{ext}$ can be either positive or negative, depending on whether the 
interior angle is less than or greater than $\pi$.\\
\noindent The following statements and formula are found in \cite{Do Carmo}, pp.265-268.\\
Let $\gamma:[0,l]\rightarrow \mathbb{R}^2$ be a positively oriented parameterization of a Jordan curve where\\
a.) $\gamma$ is one-to-one.\\
b.) Let $0=t_0,t_1, \ldots ,t_n=l$ be a partition of $[0,l]$ where $\gamma$ is $C^2$ on $(t_0,t_1)\cup \ldots \cup(t_{n-1},t_n)$ 
and $\gamma^\prime(t) \neq 0$ for all $t$ in this range.  Also let $\gamma_k$ have non-zero left and right derivatives 
at $t_1,\ldots,t_{n-1}$ with a non-zero right derivative at $t_0$ and a non-zero left derivative at $t_n$.
  Let $\theta_0\ldots \theta_{n-1}$ denote the external angles of $\gamma$ where
$\phi_i :[t_i,t_{i+1}] \rightarrow \mathbb{R}$ is a function (described in Theorem~\ref{theta})
 which describes the rotation of the tangent vector on $[t_i,t_{i+1}]$.\\
If the above is true, the following formula holds:
\begin{equation}\label{turningtangent}
\sum_{k=0}^{n-1} \left(\phi_k(t_{k+1})-\phi_k(t_k)\right)+\sum_{k=0}^{n-1} \theta_k =2\pi.
\end{equation}
\begin{figure}
\setlength{\unitlength}{1cm}
\begin{picture}(6.0,6.0)(-3.3,-5)
\rotatebox{270}{\includegraphics{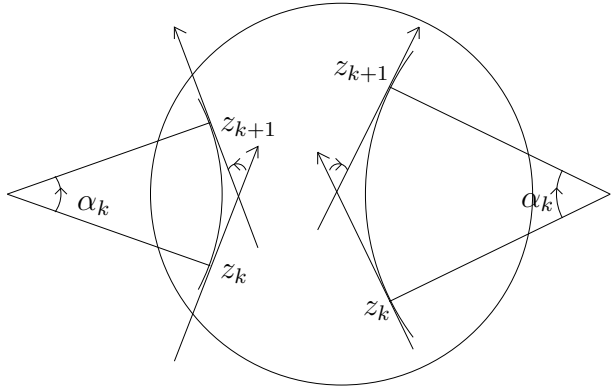}}
\end{picture}
\put(1.7,1.4){$z_k$}
\put(1.3,4.6){$z_{k+1}$}
\put(-.2,1.9){$z_k$}
\put(-.2,3.85){$z_{k+1}$}
\put(3.8,2.85){$\alpha_k$}
\put(-2.1,2.8){$\alpha_k$}
\caption{Illustration of $\phi_k(t_{k+1})-\phi_k(t_k)$ for outward (left) and inward (right) sides.}\label{inwardoutward}
\end{figure}

\noindent This theorem, though ultimately topological, states the intuitive notion that the total 
angle swept by the tangent vector, including jumps is equal to $2\pi$.\\
\noindent Let $\gamma$ be the boundary of an $n$-gon.  Refer to Figure~\ref{inwardoutward}.  In this instance, 
the geometric interpretation of $\phi_k(t_{k+1})-\phi_k(t_k)$ is $|\phi_k(t_{k+1})-\phi_k(t_k)|=\alpha_k$, where 
$\alpha_k$, $0\leq \alpha_k \leq \pi$, is the angle measure of the arc of the geodesic from 
$\gamma(t_k)$ to $\gamma(t_{k+1}).$  If the geodesic is inward, then the angular difference is negative, whereas if the 
geodesic is outward, the angular difference is positive or $0$.  This yields the following:
\begin{eqnarray}
\phi_k(t_{k+1})-\phi_k(t_k) &=& -\alpha_k \qquad\mbox{if the geodesic is inward;}\\
\phi_k(t_{k+1})-\phi_k(t_k) &=& \alpha_k \qquad\mbox{if the geodesic is outward.}
\end{eqnarray}
If, in equation~(\ref{turningtangent}), we break the first sum into the sum over 
$l$ inward sides, $\alpha_0, \ldots, \alpha_{l-1}$ and $m$ outward sides, $\beta_0, \ldots, \beta_{m-1}$
and note equation~(\ref{thetaext}) then we have
\begin{equation}
-\sum_{k=0}^{l-1} \alpha_k +\sum_{k=0}^{m-1} \beta_k +\sum_{k=0}^{n-1} \left(\pi-(\theta_{int})_k\right)=2\pi.
\end{equation}
After rearrangement and reindexing, beginning from $k=1$, we have
\begin{equation}
\pi(n-2)-\sum_{k=1}^{n} (\theta_{int})_k = \sum_{k=1}^{l} \alpha_k -\sum_{k=1}^{m} \beta_k,
\end{equation}
but the left hand side is simply equation~(\ref{classical}) for the area, so that we have proven, in an alternative fashion, 
equation~(\ref{geometric}).  Additionally, this formula is now valid for the closed unit disc.\\
\noindent  We will now derive equation~(\ref{compform}) geometrically so that it will be valid on the closed disc.
In light of criteria~(\ref{criterion1}) and~(\ref{criterion2}), all that is necessary to show is that the angle measure $\alpha_k$
 of the geodesic  between $z_k$ and $z_{k+1}$ is 
$\left|2\ {\rm Tan}^{-1} \left( {\frac{x_k y_{k+1}-y_k x_{k+1}}{1-x_k x_{k+1}-y_k y_{k+1}}} \right)\right|.$  
This is clearly the case if $z_k$ and $z_{k+1}$ lie on some ray throught the origin, i.e. when 
$x_k y_{k+1}-y_k x_{k+1}=0,$ which is equivalent to ${\rm Im}(\bar{z}_kz_{k+1})=0$, so assume that 
${\rm Im}(\bar{z}_kz_{k+1}) \neq0$.  Let $z_k,z_{k+1} \in \overline{B_1(0)^2}.$  
\begin{figure}
\setlength{\unitlength}{1cm}
\begin{picture}(6.0,6.0)(-3.3,-5)
\rotatebox{270}{\includegraphics{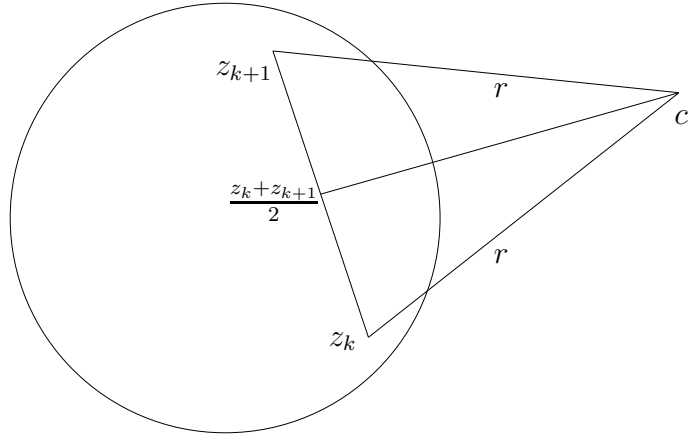}}
\end{picture}
\put(1.2,1){$z_k$}
\put(-.3,4.6){$z_{k+1}$}
\put(-.15,2.8){$\frac{z_k+z_{k+1}}{2}$}
\put(5.8,3.95){$c$}
\put(3.4,2.1){$r$}
\put(3.4,4.3){$r$}
\caption{Illustration for equation~(\ref{geometricangles}).}\label{geomangles}
\end{figure}
Refer to Figure~\ref{geomangles}.  The Euclidean line segment 
$(z_k,z_{k+1})$ has midpoint $\frac{z_k+z_{k+1}}{2}$.  The geodesic from $z_k$ to $z_{k+1}$ must be a segment of a 
circle $\mathcal{C}$ which intersects the unit disc at a right angle.  Let $c$, a complex 
number, denote the center of $\mathcal{C}$.  If $r$ is the radius of $\mathcal{C}$, then clearly $|c|^2=r^2+1.$  A 
parameterization of $\mathcal{C}$ is
\begin{equation}
z=c+\sqrt{|c|^2-1}\ e^{i\theta} \qquad\mbox{where $z$ is any point on $\mathcal{C}$ and $\theta$ is a real number.}
\end{equation}
Basic algebra shows that
\begin{equation}
\bar{z}c+z\bar{c}=1+|z|^2.
\end{equation}
In particular this is satisfied for $z=z_k$ and $z=z_{k+1}$:
\begin{eqnarray}
\bar{z}_kc+z_k\bar{c} &=& 1+|z_k|^2,\\
\bar{z}_{k+1}c+z_{k+1}\bar{c} &=& 1+|z_{k+1}|^2.
\end{eqnarray}
Solving for $c$ we obtain
\begin{equation}
c=\frac{(1+|z_k|^2)z_{k+1}-(1+|z_{k+1}|^2)z_k}{\bar{z}_kz_{k+1}-z_k\bar{z}_{k+1}}.
\end{equation}
Geometrically, we have $\alpha_k$ given by the following formula,
\begin{equation}\label{geometricangles}
\alpha_k = 2\ {\rm Tan}^{-1} \left|\frac{\frac{z_{k+1}-z_k}{2}}{c-\frac{z_k+z_{k+1}}{2}}\right|=
2\ {\rm Tan}^{-1} \left|\frac{z_{k+1}-z_k}{2c-z_k-z_{k+1}}\right|.
\end{equation}
After substitution of $c$,
\begin{equation}
\alpha_k = 2\ {\rm Tan}^{-1} \left|
\frac{(z_{k+1}-z_{k})(\bar{z}_kz_{k+1}-z_k\bar{z}_{k+1})}{2(1+|z_k|^2)z_{k+1}-2(1+|z_{k+1}|^2)z_k-
(z_k+z_{k+1})(\bar{z}_kz_{k+1}-z_k\bar{z}_{k+1})}\right|.
\end{equation}
Simplification yields
\begin{equation}
\alpha_k =2\ {\rm Tan}^{-1} \left|
\frac{(z_{k+1}-z_{k})(\bar{z}_kz_{k+1}-z_k\bar{z}_{k+1})}{(z_{k+1}-z_k)(2-z_{k+1}\bar{z}_k-\bar{z}_{k+1}z_k)}\right|.
\end{equation}
We now have the following,
\begin{equation}
\alpha_k = 2\ {\rm Tan}^{-1}\left|\frac{{\rm Im}(\bar{z}_kz_{k+1})}{{\rm Re}(1-\bar{z}_kz_{k+1})}\right|
= \left|2\ {\rm Tan}^{-1} \left( {\frac{x_k y_{k+1}-y_k x_{k+1}}{1-x_k x_{k+1}-y_k y_{k+1}}} \right)\right|.
\end{equation}
This shows that equation~(\ref{compform}) is indeed valid for the closed disc.

\chapter[MISCELLANEOUS IDENTITIES AND APPLICATIONS]{MISCELLANEOUS IDENTITIES INVOLVING AREA AND 
PERIMETER}
\setcounter{equation}{0}
\noindent In this chapter, several identities will be developed which hinge on the analytic method of Chapter II, in 
particular equations~(\ref{pureimaginary}) and~(\ref{areacomplex2}).  Among these identities is a relatively simple 
expression which couples area and perimeter, as well as two alternative formulas for area, depending upon whether 
$0$ is in the interior or exterior of the $n$-gon.

\noindent 
If $\gamma$ is the parameterization for the $n$-gon, then after calculation of equation~(\ref{perimeter}), 
the perimeter is
\begin{equation}
P = \sum_{k=1}^n \ln\left(\frac{|a_k|+1}{|a_k|-1}\right) = \sum_{k=1}^n 2\ {\rm Tanh}^{-1} \frac{1}{|a_k|}.
\end{equation}
We can rewrite this as
\begin{eqnarray}
P &=& \sum_{k=1}^n \ln\left(\frac{(|a_k|+1)^2}{|a_k|^2-1}\right)=
\sum_{k=1}^n \ln\left(\left(\frac{|a_k|+1}{|a_k|}\right)^2\left(\frac{|a_k|^2}{|a_k|^2-1}\right)\right)\nonumber\\
 &=& \sum_{k=1}^n \ln\left(\frac{|a_k|+1}{|a_k|}\right)^2+\sum_{k=1}^n \ln\left(\frac{|a_k|^2}{|a_k|^2-1}\right),
\end{eqnarray} so that we have
\begin{equation}
\sum_{k=1}^n \ln\left(\frac{|a_k|^2}{|a_k|^2-1}\right)=P+ \sum_{k=1}^n \ln\left(\frac{|a_k|}{|a_k|+1}\right)^2.
\end{equation}
When substituted into equation~(\ref{areacomplex}) this yields
\begin{equation}
A=\frac{1}{i}\sum\limits_{k=1}^n 2\log\left(1+\frac{\bar{z}_k}{a_k}\right)+\frac{P}{i}+\frac{1}{i}
\sum_{k=1}^n \ln\left(\frac{|a_k|}{|a_k|+1}\right)^2,
\end{equation}
which simplifies to
\begin{equation}\label{areaperimeter}
A+iP=\frac{2}{i}\sum_{k=1}^n \log\left(\frac{|a_k|^2+\overline{a_kz_k}}{|a_k|^2+|a_k|}\right).
\end{equation}
Surprisingly, we have a nice expression with the real part equal to area and the imaginary part equal to perimeter.\\
Equation~(\ref{compform}) rewritten as
\begin{equation}
A ={\rm Im}\sum\limits_{k=1}^n 2\log\left(1+\frac{\bar{z}_k}{a_k}\right)\,
\end{equation}
also yields two more area formulas:
\begin{eqnarray}
A &=& -4\pi+\sum\limits_{k=1}^n 2\arg\left(1+\frac{1}{a_kz_k}\right)\label{int} \quad\mbox{if 0 is in the interior of the $n$-gon} \\
A &=& \sum\limits_{k=1}^n 2\arg\left(1+\frac{1}{a_kz_k}\right)\label{ext} \qquad\mbox{if 0 is in the exterior of the $n$-gon.}
\end{eqnarray}
To prove this result we need the following claim:\\

\noindent \textit{Claim:}  Let $a \in \mathbb{C}\setminus[-1,0]$ then $\int_0^1 \frac{{d}t}{t+a}=\log\left(1+\frac{1}{a}\right)$, 
where $\log$ denotes the principal branch of the logarithm.\\

\noindent\textit{Proof of claim:}  If $a \in \mathbb{C}\setminus[-1,0]$ and $a$ is real, then clearly the claim is true.  
If $a \in \mathbb{C}\setminus[-1,0]$ and $a$ is complex, then $t+a$ is non-negative for $t \in [0,1]$, so that 
we can integrate with the principal branch of the logarithm.  We have that $\int_0^1 \frac{{d}t}{t+a}=\log(1+a)-\log(a).$  The 
claim is proven if we can show that $\arg(1+a)-\arg(a) \in (-\pi,\pi).$  Now ${\rm Im}\ (1+a)={\rm Im}\ (a)$ implies
\begin{eqnarray}
\arg(1+a) \in [0,\pi) \quad & \Longleftrightarrow & \arg(a) \in [0,\pi), \quad\mbox{and}\nonumber\\
\arg(1+a) \in (-\pi,0) \quad &\Longleftrightarrow & \arg(a) \in (-\pi,0)\nonumber.
\end{eqnarray}
In either case, it is clear that $\arg(1+a)-\arg(a) \in (-\pi,\pi).$  This concludes the proof of the claim.\\

\noindent To prove equations~(\ref{int}) and~(\ref{ext}) consider $I = \int\limits_\gamma \frac{1}{z}\mathrm {d}z$, then 
$I =\sum_{k=1}^n \int_0^1 \frac{\gamma_k^\prime(t)}{\gamma_k(t)}\mathrm {d}t$ where 
$\gamma_k(t)=\frac{t+a_kz_k}{t\bar{z}_k+a_k}$ and $t\in [0,1].$  Applying algebra
\begin{equation}
I=\sum_{k=1}^n \int_0^1 \frac{a_k(1-|z_k|^2)}{(t+a_kz_k)(t\bar{z}_k+a_k)}\mathrm {d}t
\end{equation}
and partial fractions,
\begin{equation}
I=\sum_{k=1}^n \int_0^1 \frac{1}{t+a_kz_k}\mathrm {d}t-
\sum\limits_{k=1}^n \int_0^1 \frac{\bar{z}_k}{t\bar{z}_k+a_k}\mathrm {d}t.
\end{equation}
By the claim and previous work this becomes
\begin{equation}
I=\sum_{k=1}^n \log\left(1+\frac{1}{a_kz_k}\right)-\sum\limits_{k=1}^n \log\left(1+\frac{\bar{z}_k}{a_k}\right).
\end{equation}
Taking the imaginary part of each side of the equation and rearranging terms we have
\begin{equation}
A=-2\ {\rm Im}(I)+\sum\limits_{k=1}^n 2\arg\left(1+\frac{1}{a_kz_k}\right).
\end{equation}
Now $I=2\pi i n(\gamma)$ where $n(\gamma)$ is the wrapping number of $\gamma$ about $0$, so in the first case, 
$I=2\pi i$ and in the second case $I=0$ which proves the result.\\

\noindent These identities are interesting because they perhaps suggest that sufficient algebraic 
relations exist between area and perimeter to yield a purely analytic solution to the isoperimetric problem:  Given an 
$n$-gon with fixed perimeter, maximize (if such a maximum exists) the area of the $n$-gon.  This has been proven by 
\cite{Bezdek}, but a proof of existence of a solution is not given.
By a personal communication from Dr. Roger Barnard via Dr. Petros Hadjicostas (Texas Tech University), 
we seek to show that 
\begin{equation}
A\leq \pi(n-2)-2n\;{\rm Sin}^{-1}\left(\frac{\cos(\pi/n)}{\cosh(P/2n)}\right),
\end{equation} which with minimal effort can be shown to be equivalent to
\begin{equation}
{\rm Re} \cos\left(\frac{A+2\pi}{2n}+i\frac{P}{2n}\right) \geq \cos(\pi/n).
\end{equation}
This is similar in spirit to the above identities.\\

\noindent Another use of equation~(\ref{compform}) is the proof of existence of a solution to the isoperimetric problem.
Denote $r \in B_1(0)^n$ by $r=(z_1,\ldots,z_n)$.  Denote the area equation~(\ref{compform}) by 
$A(r)$ and equation~(\ref{perimeter}) by $P(r)$.  The area and perimeter functions are certainly defined for any such 
configuration $r \in B_1(0)^n$, though the points $ z_1,\dots ,z_n$ may correspond to self intersecting 
$n$-gons or $n$-gons with negative orientation. Working on this expanded set, given a fixed perimeter, the maximum value 
of the area function exists.  This is expressed by the following theorem.
\begin{theorem}
If $P >0$ is fixed, then $A(r)$ attains its supremum over the set\\ 
$S=\{r \in B_1(0)^n : P=P(r)\}$.
\end{theorem}

\noindent \textit{Proof.}  Define 
$S_1=\{r \in B_1(0)^n : P=P(r), z_1=0\}$.  Denote 
the range of $A$ over $S$ and $S_1$ as $A(S)$ and $A(S_1)$, respectively.  Also $S_1 \subset S$ 
implies $A(S_1) \subset A(S)$.  We will now show the opposite inclusion.  Let $a \in A(S)$, then there exist 
$r_a =(z_1,\ldots ,z_n) \in S$ such that $ a=A(r_a) $.  Let $M(z)$ be an isometry of the disc such that $M(z_1)=0$.  
Define $M(r_a)=M(z_1,\ldots,z_n)$ to be $(M(z_1),\ldots,M(z_n))$.  Then $M(r_a)$ can be written as 
$M(r_a)=(0,z_2^\prime,\ldots,z_n^\prime)=r_a^\prime \in S_1.$ M is an isometry, meaning that in this case, 
$a=A(r_a)=A(M(r_a))=A(r_a^\prime) \in A(S_1)$.  This shows that $A(S) \subset A(S_1)$, that is $A(S)=A(S_1).$  Note that $A(r)$ 
is continuous on $S_1$ .\\
If we show that $S_1$ is closed, then $S_1$ is compact and $A(S_1)$ (hence $A(S)$) attains its 
supremum and we have proven the theorem.  In order to prove that $S_1$ is closed, we first need to show that $  
\overline{S_1} \subset B_1(0)^n.$  Let $r=(z_1=0,z_2,\ldots,z_n) \in S_1.$  Denote the hyperbolic distance between two 
points $z_j,z_k$ by $d(z_j,z_k)$.  Explicitly, $d(z_j,z_k)=2\ {\rm Tanh}^{-1}\left|\frac{z_k-z_j}{1-\bar{z}_jz_k}\right|.$  
Let $j$ be an integer from 1 to $n$, then by the triangle inequality and the fact that $P=\sum_{k=1}^n d(z_k,z_{k+1})$,
\begin{eqnarray}
d(0,z_j) &\leq & d(0,z_2)+\ldots +d(z_{j-1},z_j)\nonumber\\
& \leq & P.
\end{eqnarray} 
This yields $2\ {\rm Tanh}^{-1}|z_j| \leq P$, implying $|z_j| \leq \tanh P/2$, so that $z_j \in B_{\tanh P/2}(0)$.  Now 
$r \in B_{\tanh P/2}(0) ^n$, results with $S_1 \subset B_{\tanh P/2}(0) ^n$, hence
$\overline{S_1} \subset \overline{B_{\tanh P/2}(0) ^n}$. This means $\overline{S_1} \subset B_1(0)^n$ 
since for all real $P$, $\tanh P/2 \leq 1.$\\
We will now show that $S_1$ is closed.  Let $l=(z_1,\ldots,z_n)$ be a limit point of $S_1$, then there exists 
$\{r_m\}_{m=1}^\infty \subset S_1$ such that $\lim_{m\to\infty} r_m=l.$ For each $r_m$ the first coordinate is $0$, 
so that $z_1=0$, the first condition for $l \in S_1$.  Now $P(r)$ is continuous at $l$ so 
$\lim_{m\to\infty} P(r_m)=\lim_{m\to\infty} P=P = P(l)$ (the second condition for $l \in S_1$), therefore $S_1$ is closed, 
and the theorem is proven.  $\Box$

\chapter{CONCLUDING REMARKS}
\setcounter{equation}{0}

\noindent The two proofs of equation~(\ref{compform}) and its geometric interpretation equation~(\ref{geometric}) 
have respective strengths and weaknesses.  The analytic derivation (Chapter II) has the advantage of yielding some 
unexpected relationships between area and perimeter, and even (with relatively little work) alternative formulas for the 
area provided that the origin does not lie on the boundary of the $n$-gon. This proof relies on the fact that hyperbolic 
geometry is generated by a conformal metric.  The drawback to this method is that the proof is only valid if all of the 
vertices are in the interior of the unit disc, otherwise the conformal metric has a singularity which destroys the validity 
of Green's theorem.  This proof, in particular equation~(\ref{pureimaginary}), allows for some unexpected identities 
involving the perimeter and the area, and suggests a direct link to the solution of the isoperimetric problem.  One such
 link is the existence of a maximum configuration if the set of Jordan $n$-gons is expanded to include self intersections.  
Potentially, this method could also be applied to geometries generated by other conformal metrics, for instance, higher 
dimensional hyperbolic and spherical geometries by \cite{Ratcliffe}.  The geometric proof (Chapter III) relies only on the equations~(\ref{classical}),~(\ref{turningtangent}) and a 
geometric knowledge of the geodesics.  The benefit of this proof is its brevity, and it avoids the singularity of the 
conformal metric.  This proof establishes equations~(\ref{compform}) and~(\ref{geometric}) on the entire unit disc, but does 
not suggest any of the identities developed from the analytic method.    In conclusion, we see that 
apart from their inherent appeal, equations~(\ref{compform}) and~(\ref{geometric}) have ranging applicability to 
hyperbolic geometry as a whole.


\begin{thebibliography}{1000}
\bibitem{Anderson}  Anderson, J. W.  \textit{Hyperbolic Geometry.}  London, Springer. 1999.
\bibitem{Bezdek}  Bezdek, K. \textit{Ein elementarer Beweis f\"{u}r die isoperimetrische Ungleichung in der Eucklidischen 
und hyperbolischen Ebene.}  [\textit{An elementary proof for the isoperimetric inequality in the Euclidean and 
Hyperbolic planes.}]  Ann. Univ. Sci. Budapest E\"{o}tv\"{o}s Sect. Math. \textbf{27} (1984) 107-112 (1985).
\bibitem{Do Carmo}  Do Carmo, M. P.   \textit{Differential Geometry of Curves and Surfaces.}  Upper Saddle River, 
New Jersey, Prentice-Hall Inc.  1976.
\bibitem{Needham}  Needham, T. \textit{Visual Complex Analysis.} New York, New York, Oxford University Press. 1997.
\bibitem{Ratcliffe}  Ratcliffe, J. G. \textit{Foundations of Hyperbolic Manifolds.} New York, Springer-Verlag.  1994.
\bibitem{Stillwell}  Stillwell, J.  \textit{Geometry of Surfaces.}  New York, Springer-Verlag. 1992.



\end{thebibliography}
\end{document}